\documentclass{amsart}
\usepackage{amssymb,amsmath,amsthm}
\begin{document}
\title{Goldbach's Problem in Primes  with Binary \\ Expansions of a Special Form}%
\author{K.M.\'Eminyan}%
\address{Financial  University  under the Government  of the Russian Federation.  
Bauman State Technical University. Moscow.}%
\email{eminyan@mail.ru}%


\begin{abstract}
Let $\mathbb{N}_0$ be a class of natural numbers whose binary expansions contain  even numbers of ones. Goldbach's problem in numbers of class $\mathbb{N}_0$ is solved.
\end{abstract}
\maketitle
% ----------------------------------------------------------------
\begin{center}
{\bf{1. Introduction}}
\end{center}
\bigskip

Let $n=e_{0}+e_{1}2+\ldots+e_{k}2^{k}$  be a binary expansion of a natural number $n$,  $(e_{j}=0, 1)$. Let $\mathbb{N}_{0}$  be a set of natural numbers whose binary expansions have an even number of ones, $\mathbb{N}_{1}=\mathbb{N}\setminus\mathbb{N}_{0}$. Let
$$
\varepsilon(n)=\left\{
  \begin{array}{ll}
    {~~\,\,}1, & \text{if}{~~~} n\in \mathbb{N}_{0}; \\
    -1, & \text{if} {~~~}n\in \mathbb{N}_{1}.
  \end{array}
\right.
$$

In 1968,  A.O. Gelfond \cite{Gelfond} proved that numbers from the sets $\mathbb{N}_{0}$~and~$\mathbb{N}_{1}$ are regularly distributed in arithmetical progressions.

In 1991, The author got \cite{Eminyan} the asymptotical formula for the sum
$$
\sum\limits_{n\leqslant x,\,\, n\in \mathbb{N}_{0}} \tau(n)
$$
and so solved  Dirichlet divisors problem in the numbers of class
~$\mathbb{N}_{0}$.

In 2010, C. Mauduit and  J. Rivat \cite{MR} proved in particular that the densities of sets of primes of the classes $\mathbb{N}_{0}$~and~$\mathbb{N}_{1}$ are equal to each other. B. Green gave another proof of this fact   \cite{Green}. These papers are based on estimates of exponential sums of a special type, which, by the force and by methods of proofs, are variants of  estimate, derived by the author in 1991, of the integral of  modulus of a trigonometric sum of the special type \cite{Eminyan}.

In this paper the ternary Goldbach problem in prime numbers of the set $\mathbb{N}_{0}$ is solved.

The main results are contained in the following theorems.

\textbf{Theorem 1}. \textit{Let $\alpha$ be an arbitrary real number. There exists an absolute constant $\varkappa>0$ such that}
 $$
 S=\sum\limits_{n\leqslant X}\varepsilon(n)\Lambda(n)e^{2\pi i \alpha n}=O(X^{1-\varkappa}).
 $$

 \textit{The constant in sign $O$ is absolute.}

 \textbf{Theorem 2}. \textit{Let $J(N)$ be the number of representations of odd $N$ by sum of three primes, and  $J_0(N)$ be the number of representations of odd $N$ by sum of three primes from the set}~$\mathbb{N}_{0}$.

 \textit{Then the equality }
 $$
 J_0(N)=\frac{1}{8}J(N)(1+O(N^{-\varkappa}\ln N)),
 $$
\textit{holds, where $\varkappa>0$ is a constant from theorem 1.}

\bigskip
\begin{center}
{\bf{2. Auxiliary lemmas}}
\end{center}
\bigskip

\textbf{Lemma 1}. \textit{Let}
$$
\alpha=\frac{a}{q}+\frac{\theta}{q^{2}}, {~~}(a,\,q)=1,{~~}q\geqslant 1,{~~} |\theta|\leqslant 1.
$$
\textit{Then for any $\beta \in \mathbb{R}$, $U>0$, $P\geqslant 1$ we have}
$$
\sum\limits_{x=1}^P \min\left(U,\,\|\alpha x+\beta\|^{-1}\right)\leqslant 6 \left(\frac{P}{q}+1\right)(U+q\log q).
$$
Proof see in \cite[chapter 4]{AAK}.

\textbf{Lemma 2}. (A.O. Gelfond) \textit{Let $Q\in \mathbb{N}$.  The inequality}
$$
\left|\prod\limits_{r=0}^{2Q-1}\left(1-e^{2\pi i \alpha 2^{r}}\right)\right|\leqslant \frac{2}{\sqrt{3}}2^{2Q\lambda},
$$
\textit{holds, where} $\lambda=\frac{\ln 3}{\ln 4}=0,7924812\ldots$

Proof see in  \cite{Gelfond}.

\textbf{Corollary 1}. \textit{For any $\alpha \in \mathbb{R}$ the estimate}
$$
\left| \sum\limits_{n\leqslant X}\varepsilon (n)e^{2 \pi i \alpha n}\right|=O(X^{\lambda}\ln X).
$$
\textit{holds.}

Proof. Define natural number $Q$ with inequalities
$$
2^{2(Q-1)}<X+1\leqslant 2^{2Q}.
$$
Then
$$
\left| \sum\limits_{n\leqslant X}\varepsilon (n)e^{2 \pi i \alpha n}\right|=
\left| \sum\limits_{n<2^{2Q}}\varepsilon (n)e^{2 \pi i \alpha n}
\sum\limits_{n_{1}\leqslant X}\frac{1}{2^{2Q}}\sum\limits_{l=1}^{2^{2Q}}e^{2\pi i \frac{(n-n_{1})l}{2^{2Q}}}\right|\leqslant
$$
$$
\leqslant 2^{-2Q}\sum\limits_{l=1}^{2^{2Q}}\left| \sum\limits_{n<2^{2Q}}\varepsilon(n)e^{2 \pi i (\alpha+l\,2^{-2Q})n}\right|\left| \sum\limits_{n_{1} \leqslant X}e^{-2 \pi i n_{1}l 2^{-2Q}}\right|.
$$

Furthermore, it follows from from inequality
$$
\sum\limits_{n<2^{2Q}}\varepsilon(n)e^{2 \pi i (\alpha+l\,2^{-2Q})n}= \prod\limits_{r=0}^{2Q-1}\left(1-e^{2\pi i (\alpha+l 2^{-2Q})2^r}\right)
$$
and lemma 2 that
$$
\left| \sum\limits_{n\leqslant X}\varepsilon(n)e^{2\pi i \alpha n}\right|\ll X^{\lambda}2^{-2Q}\sum\limits_{l=1}^{2^{2Q}}\min(X,\,\|l 2^{-2Q}\|^{-1}).
$$

From this and Lemma 1 we have Corollary 1.

\textbf{Lemma 3}. (Gallagher). \textit{Let $S(t)$ be a  complex valued function with continuous first derivative on   $[t_0,\, t_k]$ and, $t_0<t_1<\ldots <t_{k-1}<t_k$.}

\textit{Then, assuming that $\delta=\min\limits_{0\leqslant r<k} (t_{r+1}-t_r)$,we have}

$$
\sum\limits_{r=1}^{k}|S(t_r)|\leqslant \frac{1}{\delta}\int_{t_0}^{t_k}|S(t)|\,dt+ \frac{1}{2}\int_{t_0}^{t_k}|S'(t)|\,dt.
$$
Proof see in \cite[chapter 1]{Montgomery}.

\bigskip
\begin{center}
{\bf{3. The main lemma and its corollaries}}
\end{center}
\bigskip

\textbf{Lemma 4}. \textit{Let} $Q\in \mathbb{N}$,
$$
S_{Q}(a)=\prod_{r=0}^{2Q-1}\left(1-e^{2\pi i \alpha 2^r}\right).
$$

\textit{The inequality}
$$
\int_{0}^{1}|S_Q(\alpha)|\,d\alpha\leqslant 2^{Q\theta_0},
$$
\textit{holds, where} $\theta_0=\log_2 \sqrt{2+\sqrt{2}}=0,88577\ldots$

Proof see in \cite{Eminyan}.

\textbf{Corollary 2}. \textit{The inequality}
$$
\int_{0}^{1}\left| \sum\limits_{n\leqslant X}\varepsilon (n)e^{2 \pi i \alpha n}\right|\,d\alpha \leqslant X^{\theta/2}\ln X.
$$
\textit{holds}

Proof. Let $Q\in \mathbb{N}$, $2^{2(Q-1)}<X+1\leqslant 2^{2Q}$

Then
$$
\int_{0}^{1}\left| \sum\limits_{n\leqslant X}\varepsilon (n)e^{2 \pi i \alpha n}\right|\,d\alpha \ll 2^{-2Q}\sum\limits_{l=1}^{2Q}\min (X,\,\|l 2^{-2Q}\|^{-1})\int_{0}^{1}\left|S_Q(\alpha+l 2^{-2Q}) \right|\,d\alpha.
$$

Since $S_Q(\alpha)$ is a periodic function of $t$ with period 1,
$$
\int_{0}^{1}\left|S_Q(\alpha+l 2^{-2Q}) \right|\,d\alpha=\int_{0}^{1}\left|S_Q(\alpha) \right|\,d\alpha.
$$

By lemma 1,
$$
2^{-2Q}\sum\limits_{l=1}^{2Q}\min (X,\,\|l 2^{-2Q}\|^{-1})\ll \ln X.
$$

The assertion of Corollary 2 follows from Lemma 4.

\textbf{Corollary 3}.\textit{ ѕусть $k\in \mathbb{N}$. The estimate}
$$
2^{-k}\sum\limits_{r=0}^{2^{k}-1}\left|\sum\limits_{x=0}^{2^{k}-1}\varepsilon(x)e^{2 \pi i \frac{rx}{2^{k}}} \right|\ll k 2^{\theta k}
$$
holds.

Proof.  Applying Lemma 3, putting it
$$
S(t)=\sum_{x=0}^{2^{k}-1}\varepsilon(x)e^{2\pi i t x},{~~~} t_{r}=\frac{r}{2^{k}}.
$$

Then
$$
S'(t)=2\pi i \sum_{x=0}^{2^{k}-1}x\varepsilon(x)e^{2\pi i t x},
$$
$\delta=2^{-k}$.

Then
$$
2^{-k}\sum\limits_{r=0}^{2^{k}-1}\left|\sum\limits_{x=0}^{2^{k}-1}\varepsilon(x)e^{2 \pi i \frac{rx}{2^{k}}} \right|=2^{-k}\sum\limits_{r=0}^{2^{k}-1}\left|S(t_{r})\right|\leqslant
$$
$$
\leqslant\int_{0}^{1}|S(t)|\,dt+ 2^{-k}\int_{0}^{1}|S'(t)|\,dt.
$$

Apply Abel's transform to $S'(t)$:
$$
2^{-k}\int_{0}^{1}|S'(t)|\,dt\ll 2^{-k}\int_{0}^{1}2^{k}\left| \sum\limits_{x=0}^{2^{k}-1}\varepsilon(x)e^{2 \pi i x t}\right|\,dt+
$$
$$
+2^{-k}\int_{0}^{1}\int_{0}^{2^k-1}\left(\left|\sum\limits_{x\leqslant u}\varepsilon(x)e^{2 \pi i x t}\right|\,dt\right)du\ll\int_{0}^{1}\left|\sum\limits_{x\leqslant u_{0}}\varepsilon(x)e^{2 \pi i x t}\right|dt,
$$
where $u_{0}$ is a number from the segment $[0,\,2^k-1]$ such that the last integral reaches its maximum. Now, Corollary 3 follows from Corollary 2.

\textbf{Corollary 4}. \textit{Let $k$  и $t$ be integers, $0\leqslant t \leqslant k$. Suppose that the inequality $m\in \mathbb{N}$} holds. Then
$$
\hat{\varepsilon}_{m}(r)=2^{-m}\sum_{x=0}^{2^{m}-1}\varepsilon(m)e^{-2 \pi i \frac{r x}{2^{m}}}.
$$
\textit{Let $a$ is any number of the segment $[0,\, 2^{k}-1]$.
Then}
$$
\sum_{\substack{r=0\\r\equiv a \pmod {2^t}}}^{2^k-1}|\hat{\varepsilon}_k(r)|\ll 2^{(0,5-c)(k-t)}|\hat{\varepsilon}_{t}(a)|\,k,
$$
\textit{where $c=1/2-\theta_0/2$, $\theta_0$ is a number from lemma 4.}

Proof. By definition we have
$$
\sum\limits_{{r=0}\atop{r\equiv a \pmod{2^t}}}^{2^k-1}|\hat{\varepsilon}_k(r)|=2^{-k}\sum\limits_{r\pmod{2^{k-t}}}\left|\sum\limits_{x=0}^{2^{k}-1}\varepsilon(x)e^{2\pi i \frac{a+2^{t}r}{2^{k}}x}\right|=
$$
$$
=2^{-k}\sum\limits_{r\pmod{2^{k-t}}}\left|\sum\limits_{x=0}^{2^{k-t}-1} \sum\limits_{y=0}^{2^{t}-1}\varepsilon(x+2^{k-t}y)e^{2\pi i \frac{a+2^{t}r}{2^{k}}(x+2^{k-t}y)}\right|.
$$

Since $\varepsilon(x+2^{k-t}y)=\varepsilon(x)\varepsilon(y)$, we have
$$
\sum\limits_{{r=0}\atop{r\equiv a (\mod 2^t)}}^{2^k-1}|\hat{\varepsilon}_k(r)|=2^{-(k-t)}\sum\limits_{r(\mod 2^{k-t})}\left|\sum\limits_{x=0}^{2^{k-t}-1}\varepsilon(x)e^{2\pi i \left({\frac{r}{2^{k-t}}+ \frac{a}{2^{k}}}\right)x}\right|\left|2^{-t}\sum\limits_{y=0}^{2^{t}-1}\varepsilon(y)e^{2\pi i \frac{ay}{t}}\right|=
$$
$$
=2^{-(k-t)}\sum\limits_{r(\mod 2^{k-t})}\left|\sum_{x=0}^{2^{k-t}-1}\varepsilon(x)e^{2\pi i \left(\frac{r}{2^{k-t}}+\frac{a}{2^{k}}\right)x}\right||\varepsilon_t (a)|.
$$

The sum on the right side of this inequality is estimated in the same way as a similar amount of Corollary 3.

\bigskip
\begin{center}
{\bf{4. Proof of Theorem 1}}
\end{center}
\bigskip

Using Vaughan's identity (see, eg \cite[chapter 3, problem 9]{AAK}), with $u=X^{0,1}$:
$$
S=\sum_{n\leqslant X}\Lambda(n)\varepsilon(n)=W_{1}-W_{2}-W_{3}+O(u\ln u),
$$
where
$$
W_{1}=\sum\limits_{d\leqslant u}\mu(d)\sum\limits_{n\leqslant Xd^{-1}} \varepsilon(d n)e^{2\pi i \alpha dn}\ln n,
$$
$$
W_{2}=\sum\limits_{d\leqslant u}\mu(d)\sum\limits_{n\leqslant u} \Lambda(n) \sum\limits_{d n r \leqslant X}\varepsilon(d n r)e^{2\pi i \alpha d n r},
$$
$$
W_{3}=\sum\limits_{u<m \leqslant Xu^{-1}}a_{m}\sum\limits_{u<n\leqslant Xm^{-1}} \Lambda(n)  \varepsilon(m n)e^{2\pi i \alpha m n},
$$
$$
a_{m}=\sum\limits_{d|m, d\leqslant u}\mu(d).
$$

Sums $W_{1}$  and $W_{2}$ are estimated in the same way. Estimate  $W_{1}$.

Fix $d\leqslant u$. Apply to the inner sum, which we denote $S_{1}(d)$, the Abel transform, we obtain:
$$
|S_{1}(d)|\ll \left|\sum\limits_{d n\leqslant u_{0}}\varepsilon(d n) e^{2\pi i \alpha d n} \right|\log X,
$$
where $u_{0}$ is a number not exceeding $X$.

Furthermore,
$$
\sum_{dn \leqslant u_{0}} \varepsilon(d n) e^{2\pi i \alpha d n}= \sum_{m \leqslant u_{0}} \varepsilon(m) e^{2\pi i \alpha m}\,\frac{1}{d} \sum\limits_{b=0}^{d-1} e^{2\pi i \frac{b m}{d} },
$$
$$
\left|\sum\limits_{dn \leqslant u_{0}}\varepsilon (d n)e^{2 \pi i \alpha d n}\right|\leqslant \frac{1}{d}\sum\limits_{b=0}^{d-1}\left|\sum\limits_{m\leqslant u_{0}}\varepsilon(m)e^{2 \pi i \left(\alpha+\frac{b}{d}\right)m}\right|.
$$

The sum over $m$  estimate by Corollary 1:
$$
\left|\sum\limits_{m\leqslant u_{0}}\varepsilon(m)e^{2 \pi i \left(\alpha+\frac{b}{d}\right)m}\right|\ll X^\lambda\ln X,
$$
where $\lambda=0,792...$

 Thus, for any $d\leqslant u$ we have
 $$
 |S_{1}(d)|\ll X^{\lambda}\ln X,
 $$
therefore,
$$
|W_{1}|\ll u X^\lambda \ln X.
$$

Similarly, we arrive at the estimate
$$
|W_{2}|\ll u^{2} X^\lambda\ln X.
$$

The parameter $u$ is chosen so that
$$
|W_{1}|\ll  X^{1-\varkappa_{1}} {~~}\text{и}{~~}|W_{2}|\ll  X^{1-\varkappa_{1}},
$$
where $\varkappa_{1}>0$ is an absolute constant.

We now estimate $W_{3} $. Divide the interval of summation over $m$ in $O(\ln X)$ intervals of the form $\left(\frac{M}{2},\,M\right]$, where $u<\frac{M}{2}<M\leqslant \frac{X}{u}$; one of these intervals may be incomplete. So, we have:
$$
W_{3}=\sum\limits_{M}^{\ln X} W_{3}(M),
$$
where
$$
W_{3}(M)=\sum\limits_{M/2<m\leqslant M_1}a_m \sum\limits_{u<n\leqslant X/m}\Lambda (n)\varepsilon(mn) e^{2 \pi i \alpha m n},
$$
where $M/2<M_1\leqslant M$.

Furthermore,
$$
W_{3}(M)=\sum\limits_{M/2<m\leqslant M_1}a_m \sum\limits_{u<n\leqslant X/M_1}\Lambda (n)\varepsilon(mn) e^{2 \pi i \alpha m n}+$$
$$+\sum\limits_{M/2<m\leqslant M_1}a_m \sum\limits_{X/M_1<n\leqslant X/m}\Lambda (n)\varepsilon(mn) e^{2 \pi i \alpha m n}.
$$

Splitting the interval of summation over $n$ in $O(\log X)$ intervals of the form
$\left(\frac{N}{2},\, N_1\right]$, where $\frac{N}{2}<N_1 \leqslant N $, $u<N\leqslant\frac{X}{M_1}$, we arrive at the inequality
$$
|W_3|\ll |W_{3}(M,\,N)|\ln^{2}X,
$$
where
$$
W_{3}(M,\,N)=\sum\limits_{M/2<m\leqslant M_1}a_m \sum\limits_{N/2<n\leqslant N_{1}}\Lambda (n)\varepsilon(mn) e^{2 \pi i \alpha m n},
$$
where
$u<M/2<M_1\leqslant M\leqslant Xu^{-1}$,  $u<N/2<N_1 \leqslant N \leqslant Xu^{-1}$; it may be that $N_1=X/m $.

Without loss of generality, we assume that  $M \leqslant N$.

Using the fact that $|a_m|\leqslant \tau(m)\ll m^\varepsilon$  we apply the Cauchy inequality:
$$
|W_{3}(M,\,N)|^2\ll M^{1+\varepsilon}\sum\limits_{M/2<m\leqslant M_1}\left|\sum\limits_{N/2<n\leqslant N_{1}}\Lambda (n)\varepsilon(mn) e^{2 \pi i \alpha m n}\right|^2.
$$
Let $H=[X^\rho]$, where $0<\rho<10^{-3}$ - a small parameter to be chosen later. We apply van der Corput' inequality (see, eg,\cite{Green}, \cite[глава 1]{AAK}):
$$
\left|\sum\limits_{N/2<n\leqslant N_{1}}\Lambda (n)\varepsilon(mn) e^{2 \pi i \alpha m n}\right|^2\ll \frac{N}{H}\sum\limits_{|h|\leqslant H} \left(1-\frac{|h|}{H}\right)\times $$

$$
\times\sum\limits_{{N/2<n\leqslant N_{1}},\atop{N/2<n+h\leqslant N_{1}}} \Lambda(n)\Lambda(n+h)\varepsilon(m n)\varepsilon(m n+m h) e^{2 \pi i \alpha m n}e^{-2 \pi i \alpha m (n+h)}.
$$

The contribution of $h=0$ is estimated as $O\left(\frac{N^{2}}{H}\right)$, so
$$
|W_3(M,\,N)|^2\ll \frac{X^{1+\varepsilon}}{H}\sum\limits_{h=1}^{H}\sum\limits_{N/2<n\leqslant N_{1}}\left|\sum\limits_{M/2<m\leqslant M_{1}}\varepsilon(m n)\varepsilon(m n+m h) e^{-2 \pi i \alpha m (n+h)}\right|+\frac{X^{2+\varepsilon}}{H}.
$$

Fix $h\in [1,\,H]$. Now it is sufficient to prove that
$$
W_{4}(M,\, N)=\sum\limits_{N/2<n\leqslant N_{}}\left|\sum\limits_{M/2<m\leqslant M_{1}}\varepsilon(m n)\varepsilon(m n+m h)e^{-2 \pi i \alpha m (n+h)}\right|\ll X^{1-\varkappa}.
$$

Choose a positive integer $k$ from the inequalities
$$
2^{k-1}<MX^{2\rho}\leqslant 2^k.
$$

Introduce the symbol $\varepsilon_{k}(n)$:
$$
\varepsilon_{k}(n)=\left\{
  \begin{array}{ll}
    {~~\,\,}\,1, & \hbox{\text{if the sum of the first $k$ binary digits of $n$ is even};} \\
    -1, & \hbox{\text{--- otherwise}.}
  \end{array}
\right.
$$

Prove that
$$
|W_{4}(M,\, N)|\ll \sum\limits_{N/2<n\leqslant N_{}}\left|\sum\limits_{M/2<m\leqslant M_{1}}\varepsilon_{k}(m n)\varepsilon_{k}(m n+m h)e^{-2 \pi i \alpha m (n+h)}\right|+X^{1-\rho+\varepsilon}.
$$

Divide $mn$ by $2^{k}$ with remainder: $mn=2^k q +r$, $0\leqslant r<2^k$. If $r<2^k -2MH$, then $mn+mh=2^k q+r+mh$, $0<r+mh<2^k$.

For such $mn$, we have
$$
\varepsilon(mn)\varepsilon(mn+mh)=\varepsilon(r)\varepsilon(r+mh)=\varepsilon_{k}(mn)\varepsilon_{k}(mn+mh).
$$

The number of pairs $(m,\,n)$ such that the remainder $r$ lies  between $2^k -2MH$ and $2^k -1$ is $O(X2^{-k} MH X^\varepsilon)=O(X^{1-\rho+\varepsilon})$.

Now we  estimate
$$
|W_{5}(M,\,N)|=\sum\limits_{N/2<n\leqslant N_{}}\left|\sum\limits_{M/2<m\leqslant M_{1}}\varepsilon_{k}(m n)\varepsilon_{k}(m n+m h) e^{-2 \pi i \alpha h m}\right|.
$$

Introduce the discrete Fourier transform for the character $\varepsilon_{k}(r)$:
$$
\hat{\varepsilon}_{k}(r)=2^{-k}\sum\limits_{l=0}^{2^k -1}\varepsilon_k (l)e^{-2 \pi i \frac{rl}{2^k}}.
$$

From this definition it follows that
$$
\varepsilon_{k}(mn)=\sum\limits_{r=0}^{2^k -1}\hat{\varepsilon}_k (r)\exp\left\{\frac{2 \pi i r m n}{2^{k}}\right \},\
\varepsilon_{k}(mn+mh)=\sum\limits_{s=0}^{2^k -1}\hat{\varepsilon}_k (s)\exp\left\{\frac{2 \pi i s (m n+m h)}{2^{k}}\right \}.
$$

Summing  linear sums over $m$, we obtain:
$$
|W_{5}(M,\,N)|\leqslant \sum\limits_{r=0}^{2^k -1}\sum\limits_{s=0}^{2^k -1}|\hat{\varepsilon}_{k}(r)||\hat{\varepsilon}_{k}(s)|\sum_{N/2<n\leqslant N}\min\left(M,\, \left\|\frac{r+s}{2^{k}}n+\frac{hs}{2^{k}}-h\alpha\right\|^{-1}\right).
$$

From now on we will assume that $\rho=\frac{c}{200}$, where $c=\frac{1-\theta_0}{2}$, $\theta_0$ is a constant from lemma~4.

Let $t$ -- non-negative integer such that $2^t \|(r+s)$.

Suppose first $0\leqslant t \leqslant k -\frac{2\rho}{c}\log_2 X$. Note that from inequality
$2^k >M\geqslant X^{1/10}$, it follows that
 с $k>\frac{1}{10}\log_2 X$; this and $\frac{2\rho}{c}=\frac{1}{100}$ implies the inequality $k-\frac{2\rho}{c}\log_2 X\geqslant \frac{4}{5}k$.

We apply Lemma 1  with $q=2^{k-t}$, $\alpha=\frac{r+s}{2^{k}}$, $\beta=\frac{sh}{2^{k}}-\alpha h$ to the sum
$$
\sum\limits_{N/2<n\leqslant N}\min\left(M,\, \left\| \frac{r+s}{2^{k}}n+\frac{sh}{2^{k}}-\alpha h   \right\|^{-1}\right):
$$
$$
\sum\limits_{N/2<n\leqslant N}\min\left(M,\, \left\| \frac{r+s}{2^{k}}n+\frac{sh}{2^{k}}-\alpha h   \right\|^{-1}\right)\ll \left(\frac{N}{2^{k-t}}+1\right)\left(M+k2^{k-t}\right).
$$

Simplify the right-hand side of this inequality.

We have:
$$
\frac{N}{2^{k-t}}+1<\frac{NX^{2\rho}}{2^{k-t}}+1\leqslant 2 \frac{NX^{2\rho}}{2^{k-t}};
$$
$$
M+k2^{k-t}<MX^{2\rho}+k2^{k-t}<2MX^{3\rho};
$$
and so,
$$
\left(\frac{N}{2^{k-t}}+1\right)(M+k2^{k-t})<4MNX^{5\rho}2^{-(k-t)}\leqslant 5X^{1+5\rho}2^{-(k-t)},
$$
$$
\sum\limits_{N/2<n\leqslant N}\min\left(M,\, \left\| \frac{r+s}{2^{k}}n+\frac{sh}{2^{k}}-\alpha h   \right\|^{-1}\right)\ll X^{1+5\rho}2^{-(k-t)}.
$$

Now estimate the sum
$$\mathop{\sum\limits_{r=0}^{2^{k}-1}
\sum\limits_{s=0}^{2^{k}-1}}_{2^t\|(r+s)}|\hat{\varepsilon}_{k}(r)||\hat{\varepsilon}_{k}(s)|\leqslant \sum\limits_{a=0}^{2^{t}-1} \sum_{\substack{r=0\\r\equiv a\pmod 2^{k}}}^{2^{k}-1}|\hat{\varepsilon}_{k}(r)| \sum\limits_{\substack{s=0\\s\equiv -a\pmod 2^{k}}}^{2^{k}-1}|\hat{\varepsilon}_{k}(s)|.
$$

Use the corollary 4:
$$
\mathop{\sum\limits_{r=0}^{2^{k}-1}
\sum\limits_{s=0}^{2^{k}-1}}_{2^t\|(r+s)}|\hat{\varepsilon}_{k}(r)||\hat{\varepsilon}_{k}(s)|\ll 2^{(1-2c)(k-t)}\log^{2}X\sum\limits_{a=0}^{2^{k}-1}|\hat{\varepsilon}_{k}(a)|^{2}= 2^{(1-2c)(k-t)}\ln^{2}X.
$$

In this case, the estimate is achieved
$$
|W_{5}(M,\,N)|\ll X^{1+6\rho}2^{-2c(k-t)}\ln^{2}X.
$$
Recall that
$$
k-t\geqslant \frac{4}{5}k,{~~~}2^{k}>M\geqslant X^{1/10},
$$
it follows that
$$
|W_{5}(M,\,N)|\ll X^{1+6\rho-\frac{8}{50}c};
$$
finally, from the inequality $\rho\leqslant \frac{c}{100}$ we have $6\rho-\frac{8}{50}c<-\rho$,
$$
|W_{5}(M,\,N)|\ll X^{1-\rho}.
$$

It remains to consider the case
$$
k-\frac{2\rho}{c}\log_{2}X<t\leqslant k.
$$

We have:
$$
|W_{5}(M,\,N)|\leqslant \mathop{\sum\limits_{r=0}^{2^{k}-1}
\sum\limits_{s=0}^{2^{k}-1}}_{2^t\|(r+s)}|\hat{\varepsilon}_{k}(r)||\hat{\varepsilon}_{k}(s)| \min\left(M,\, \left\| \frac{r+s}{2^{k}}n+\frac{sh}{2^{k}}-\alpha h   \right\|^{-1}\right)=
$$
$$
=\sum\limits_{N/2<n\leqslant N}\sum\limits_{s_{2}=0}^{2^{k-t}-1}\sum\limits_{r_{2}=0}^{2^{k-t}-1} \mathop{\sum\limits_{s_{1}=0}^{2^{t}-1}\sum\limits_{r_{1}=0}^{2^{k}-1}}_{r_1+s_1\equiv 0\pmod{2^t}} |\hat{\varepsilon}_{k}(r_{1}+2^{t}r_{2})||\hat{\varepsilon}_{k}(s_{1}+2^{t}s_{2})|\times
$$
$$
\times\min\left(M,\, \left\| \frac{r_{1}+s_{1}+2^{t}(r_{2}+s_{2})}{2^{k}}n+\frac{(s_{1}+2^{t}s_{2})h}{2^{k}}-\alpha h   \right\|^{-1}\right).
$$

It follows from the inequalities $0\leqslant s_{1}; r_{1}< 2^t$ and congruence  $r_{1}+s_{1}\equiv 0\pmod {2^{t}}$ that either $r_{1}=s_{1}=0$ or $r_{1}+s_{1}=2^{t}$.

From this and Lemma 2 we get
$$
|W_{5}(M,\,N)|\ll 2^{-2k(1-\lambda)}\sum\limits_{N/2<n\leqslant N}\sum\limits_{s_{2}=0}^{2^{k-t}-1}
\sum\limits_{r_{2}=0}^{2^{k-t}-1} \sum\limits_{s_{1}=0}^{2^{k}-1}\min\left(M,\, \left\| \frac{hs_1}{2^k}+\beta   \right\|^{-1}\right)+
$$
$$
+X2^{-2k(1-\lambda)}2^{2(k-t)},\eqno(1)
$$
where $\beta=\frac{1+r_2+s_2}{2^{k-t}}n+\frac{hs_2}{2^{k-t}}-\alpha h$.

Let $\frac{hs_1}{2^k}=\frac{h_1s_1}{2^{k_1}}$, where $(h_1s_1,2)=1$. From Lemma 1 and inequalities
$$
2^{-k_1}\leqslant 2^{-k}X^\rho,\ \ M+2^{k_1}k_1\ll MX^{3\rho}
$$
it follows that
$$
\sum\limits_{s_{1}=0}^{2^{k}-1}\min\left(M,\, \left\| \frac{hs_1}{2^k}+\beta   \right\|^{-1}\right)\ll MX^{5\rho}.
$$

Substituting this inequality in (1):
$$
|W_{5}(M,\,N)|\ll X^{1+5\rho}2^{k-t}2^{-2k(1-\lambda)}.
$$

Now use the fact that
$$
2^{k-t}\leqslant X^{4\rho/c},\ \ \frac{\rho}{c}\leqslant \frac{1}{200},\ \ c<0.06,\ \ 2^k>M\geqslant X^{0.1}.
$$

We got:
$$
|W_{5}(M,\,N)|\ll X^{1-0.01}.
$$

Theorem 1 is proved.

\bigskip
\begin{center}
{\bf{5. Proof of Theorem 2}}
\end{center}
\bigskip

Define sums $S(\alpha)$ and $S_0(\alpha)$:
$$
S(\alpha)=\sum_{p\leqslant N}e^{2\pi i\alpha p},\ \ S_0(\alpha)=\sum_{p\leqslant N}\varepsilon(p)e^{2\pi i\alpha p}.
$$

Then
$$
J_0(N)=\frac{1}{8}\int_0^1(S(\alpha)+S_0(\alpha))^3e^{-2\pi i\alpha N }d\alpha.
$$

Expanding the brackets and using Theorem 1 and Cauchy's inequality, we obtain
$$
J_0(N)=\frac{1}{8}\int_0^1S^3(\alpha)e^{-2\pi i \alpha N }d\alpha+O(\pi(N) N^{1-\varkappa}).
$$

Since
$$
J(N)=\int_0^1S^3(\alpha)e^{-2\pi i N \alpha}d\alpha,\ \ J(N)\gg N^2(\ln N)^{-3}
$$
(for sufficiently large odd $N$), theorem 2 is proved.

\end{document}